# A unifying Lyapunov-based framework for the event-triggered control of nonlinear systems

Romain Postoyan, Adolfo Anta, Dragan Nešić, Paulo Tabuada

*Abstract*— We present a prescriptive framework for the event-triggered control of nonlinear systems. Rather than closing the loop periodically, as traditionally done in digital control, in event-triggered implementations the loop is closed according to a state-dependent criterion. Event-triggered control is especially well suited for embedded systems and networked control systems since it reduces the amount of resources needed for control such as communication bandwidth. By modeling the event-triggered implementations as hybrid systems, we provide Lyapunov-based conditions to guarantee the stability of the resulting closed-loop system and explain how they can be utilized to synthesize event-triggering rules. We illustrate the generality of the approach by showing how it encompasses several existing event-triggering policies and by developing new strategies which further reduce the resources needed for control.

## I. INTRODUCTION

The implementation of controllers on shared digital platforms offers a number of advantages in terms of cost, ease of maintenance and flexibility compared to classical dedicated control structures. However, it also poses several implementation problems, in particular we need to know when the control loop has to be closed to ensure stability and performance. In traditional setups, this is done periodically, independently of the current state of the plant. Although this approach is appealing from the analysis and implementation point of view, it often leads to unnecessary resource usage (e.g. communication bandwidth, computation time). An alternative implementation, known as event-triggered control, consists in closing the loop according to a rule that depends on the current state of the plant. A number of works addressed this topic, e.g. [2], [4], [6], [8], [15], [19], [20]. In [19], a simple strategy is proposed for nonlinear systems. The idea is the following. Assuming the continuous-time closed-loop system is input-to-state stable (ISS) with respect to measurement errors, a triggering condition is derived to guarantee that the Lyapunov function $V$ for the continuous-time system always decreases at a given rate when control tasks are executed at discrete time instants. It is shown that there does exist a constant minimal time interval between executions that reinforces the idea that event-triggered control is expected to generate larger inter-event intervals than periodic rules. This translates into a lower usage of the communication bandwidth and the computational resources. Noting that the monotonic decrease of $V$ is not necessary to guarantee asymptotic stability for the obtained hybrid systems, a triggering rule is developed in [20] to ensure that $V$ appropriately decreases at each transmission instant. This method was shown to potentially exhibit larger inter-event intervals compared to [19].

In this paper, we present a prescriptive framework for the event-triggered control of nonlinear systems. We model the problem as a hybrid system using the formalism of [7], as in [6]. We start by identifying the key features of the strategy in [19] in terms of a hybrid Lyapunov function and use them to introduce the main idea of our approach. The proposed framework relies on Lyapunov-based conditions that can be used to synthesize event-triggering rules to guarantee asymptotic stability properties. Our approach encompasses the strategies in [19], [20] for which we propose new stability analyses and relax some of the required conditions. We also develop a family of new triggering rules inspired by [20] which allows us to trade performance for longer inter-event times. In the companion paper [17], we show how this framework can be applied to distributed networked control systems subject to scheduling.

The paper is organized as follows. Notation and preliminary definitions are presented in Section II. The problem is stated in Section III where we show how the work in [19] can be analysed using a hybrid Lyapunov function. In Section IV, we develop our framework and provide guidelines on how to use it. Afterwards, we show how the work in [20] can be interpreted using our approach in Section V and propose new triggering rules. An illustrative example is provided in Section VI. The proofs are provided in the Appendix.

## II. NOTATION AND DEFINITIONS

Let $\mathbb{R} = (-\infty, \infty)$, $\mathbb{R}_{\geq 0} = [0, \infty)$, $\mathbb{R}_{>0} = (0, \infty)$, $\mathbb{Z}_{\geq 0} = \{0, 1, 2, \ldots\}$, $\mathbb{Z}_{>0} = \{1, 2, \ldots\}$. A function $\gamma : \mathbb{R}_{\geq 0} \to \mathbb{R}_{\geq 0}$ is of class $\mathcal{K}$ if it is continuous, zero at zero and strictly increasing, and it is of class $\mathcal{K}_\infty$ if in addition it is unbounded. A continuous function $\gamma : \mathbb{R}_{\geq 0}^2 \longrightarrow \mathbb{R}_{\geq 0}$ is of class $\mathcal{KL}$ if for each $t \in \mathbb{R}_{\geq 0}$, $\gamma(\cdot, t)$ is of class $\mathcal{K}$, and, for each $s \in \mathbb{R}_{>0}$, $\gamma(s, \cdot)$ is decreasing to zero.

This technical report is an extended version of [16]. This work was supported by the Australian Research Council under the Future Fellowship and Discovery Grants Schemes, the NSF awards 0820061 and 0834771 and the Alexander von Humboldt Foundation

R. Postoyan is with the Centre de Recherche en Automatique de Nancy, UMR 7039, Nancy-Université, CNRS, France romain.postoyan@cran.uhp-nancy.fr

A. Anta is with the Technische Universität Berlin & Max Planck Institute für Dynamik komplexer technischer Systeme, Germany anta@control.tu-berlin.de

D. Nešić is with the Department of Electrical and Electronic Engineering, the University of Melbourne, Parkville, VIC 3010, Australia d.nesic@ee.unimelb.edu.au

P. Tabuada is with the Department of Electrical Engineering, University of California at Los Angeles, Los Angeles, CA 90095-1594, USA tabuada@ee.ucla.edu

Additionally, a function $\beta : \mathbb{R}^3_{\geq 0} \to \mathbb{R}_{\geq 0}$ is of class $\mathcal{KLL}$, if $\beta(\cdot, \cdot, t) \in \mathcal{KL}$ and $\beta(\cdot, t, \cdot) \in \mathcal{KL}$ for any $t \in \mathbb{R}_{\geq 0}$. For $(x, y) \in \mathbb{R}^{n+m}$, the notation $(x, y)$ stands for $[x^{\mathrm{T}}, y^{\mathrm{T}}]^{\mathrm{T}}$. The distance of a vector $x$ to a set $\mathcal{A} \subset \mathbb{R}^n$ is denoted by $|x|_{\mathcal{A}} = \inf\{|x - y| : y \in \mathcal{A}\}$.

We will consider locally Lipschitz Lyapunov functions (that are not necessarily differentiable everywhere), therefore we will use the Clarke derivative which is defined as follows for $V : \mathbb{R}^n \to \mathbb{R}_{\geq 0}$: $V^\circ(x; v) := \lim\sup_{\substack{h \to 0^+ \\ y \to x}} \frac{V(y+hv) - V(y)}{h}$, that corresponds to the usual derivative when $V$ is continuously differentiable. We define the generalized gradient of $f : \mathbb{R}^n \to \mathbb{R}$ at $x$ as: $\partial f(x) := \{\zeta \in \mathbb{R}^n : f^\circ(x; v) \geq \langle \zeta, v \rangle \ \forall v \in \mathbb{R}^n\}$, that matches the classical notion of gradient when $f$ is differentiable.

We will write hybrid systems using the models proposed in [7], that are of the form:
$$\dot{x} = f(x) \quad x \in C, \qquad x^+ = g(x) \quad x \in D, \qquad (1)$$
where $x \in \mathbb{R}^n$ is the state and $C, D \subset \mathbb{R}^n$ are respectively the flow and the jump sets. Hence, any hybrid system is defined by a tuple $(C, D, f, g)$. The solutions of (1) are defined on so-called hybrid time-domains. A set $E \subset \mathbb{R}_{\geq 0} \times \mathbb{Z}_{\geq 0}$ is called a compact hybrid time domain if $E = \bigcup_{j \in \{0, \ldots, J\}} ([t_j, t_{j+1}], j)$ for some finite sequence of times $0 = t_0 \leq t_1 \leq \ldots \leq t_J$. $E$ is a hybrid time domain if for all $(T, J) \in E$, $E \cap ([0, T] \times \{0, 1, \ldots, J\})$ is a compact hybrid time domain. A hybrid signal is a function defined on a hybrid time domain. A hybrid arc is a function $\phi$ defined on a hybrid time domain $\text{dom}\,\phi$, and such that $\phi(\cdot, j)$ is locally absolutely continuous for each $j$. A hybrid arc $\phi : \text{dom}\,\phi \to \mathbb{R}^n$ is a solution to (1) if:
(1) For all $j \in \mathbb{Z}_{\geq 0}$ and almost all $t \in \mathbb{R}_{\geq 0}$ such that $(t, j) \in \text{dom}\,\phi$ we have: $\phi(t, j) \in C$, $\dot{\phi}(t, j) = f(\phi(t, j))$.
(2) For $(t, j) \in \text{dom}\,\phi$ such that $(t, j+1) \in \text{dom}\,\phi$, we have $\phi(t, j) \in D$, $\phi(t, j+1) = g(\phi(t, j))$.

Assuming $f$ and $g$ are continuous and $C, D$ are closed, system (1) possesses solutions that may be non-unique, see [7]. We are interested in the following stability definition.

**Definition 1.** *The closed set $\mathcal{A} \subset \mathbb{R}^n$ is **semiglobally asymptotically stable (S-GAS)** for system (1) if for any $\Delta \in \mathbb{R}_{>0}$ there exists $\beta_\Delta \in \mathcal{KL}$ such that for any solution $\phi$ to (1) with $|\phi(0,0)|_{\mathcal{A}} \leq \Delta$: $|\phi(t,j)|_{\mathcal{A}} \leq \beta_\Delta(|\phi(0,0)|_{\mathcal{A}}, t, j)$ for all $(t, j) \in \text{dom}\,\phi$.*

**Remark 1.** *It is shown for continuous-time systems in Proposition 3.4 in [1] that a closed set $\mathcal{A}$ is S-GAS if and only if it is globally asymptotically stable. We note that the stability bound in Definition 1 does not imply forward completeness and this is referred to as pre-asymptotic stability according to [5].*

We will show that two successive jumps (that will correspond to data transmissions in our study) are always separated by a certain uniform amount of time as long as the solution is not in the stable set $\mathcal{A}$.

**Definition 2.** *For any forward invariant set[1] $\mathcal{A} \subset \mathbb{R}^n$ for system (1), we say that **solutions to (1) have a semiglobal dwell time on** $\mathbb{R}^n \backslash \mathcal{A}$ if, for any $\Delta \in \mathbb{R}_{>0}$, there exists $\tau(\Delta) \in \mathbb{R}_{>0}$ such that for any solution $\phi$ to (1) with $|\phi(0,0)|_{\mathcal{A}} \leq \Delta$, $j \in \mathbb{Z}_{\geq 0}$:*
$$\begin{aligned}\left(\phi(t, j) \notin \mathcal{A} \quad \forall t \in [t_j, t_{j+1}]\right) \\ \Rightarrow (t_{j+1} - t_j \geq \tau(\Delta)),\end{aligned} \qquad (2)$$
*where $\text{dom}\,\phi = \bigcup ([t_j, t_{j+1}], j)$.*

## III. PROBLEM STATEMENT

### A. System models

Consider the following plant:
$$\dot{x}_P = f_P(x_P, u), \qquad (3)$$
where $x_P \in \mathbb{R}^{n_P}$ is the plant state, $u \in \mathbb{R}^{n_u}$ the control input for which a stabilizing dynamic state-feedback controller is designed:
$$\dot{x}_C = f_C(x_C, x_P), \qquad u = g_C(x_C, x_P), \qquad (4)$$
where $x_C \in \mathbb{R}^{n_C}$ is the controller state. On digital platforms, transmission between the sensors, the controller and the actuators only occur at some transmission instants $t_j$, $j \in \mathbb{Z}_{>0}$. The problem can then be modeled as follows:
$$\left.\begin{aligned}\dot{x}_P &= f_P(x_P, \hat{u}) & \forall t \in [t_{j-1}, t_j] \\ \dot{x}_C &= f_C(x_C, \hat{x}_P) & \forall t \in [t_{j-1}, t_j] \\ u &= g_C(x_C, \hat{x}_P) & \\ \dot{\hat{x}}_P &= \hat{f}_P(x_P, x_C, \hat{x}_P, \hat{u}) & \forall t \in [t_{j-1}, t_j] \\ \dot{\hat{u}} &= \hat{f}_C(x_P, x_C, \hat{x}_P, \hat{u}) & \forall t \in [t_{j-1}, t_j] \\ \hat{x}_P(t_j^+) &= x_P(t_j) & \\ \hat{u}(t_j^+) &= u(t_j) &\end{aligned}\right\} \quad (5)$$
where $\hat{x}_P$ and $\hat{u}$ denote the variables respectively generated from the most recently transmitted plant state and control input. They are usually kept constant between two transmission instants i.e. $\hat{x}_P(t) = x_P(t_{j-1})$ and $\hat{u}(t) = u(t_{j-1})$ for $t \in [t_{j-1}, t_j]$ that corresponds to $\hat{f}_P = 0$ and $\hat{f}_C = 0$. However, other implementations are possible. At each transmission instant, the controller receives $x_P(t_j)$, updates $\hat{x}_P(t_j^+) = x_P(t_j)$, sends the control input $u(t_j)$ and the actuators update $\hat{u}(t_j^+) = u(t_j)$. We suppose that this process occurs in a synchronized manner and leave the study of the effects of the eventual induced delays for future work.

Traditionally, the sequence of $t_j$, $j \in \mathbb{Z}_{>0}$, is periodic, i.e. $t_j - t_{j-1} = T$ where $T \in \mathbb{R}_{>0}$. The stability of system (5) is then guaranteed by selecting $T$ sufficiently small, see [3], [10], [13] to mention a few. In this study, we abandon this paradigm and implicitly define the transmission instants by a rule based on the states of system (5). Rewriting the problem using the hybrid formalism in [7], similar to Section II.C in [6], we obtain:
$$\left.\begin{aligned}\dot{x} &= f(x, e) \\ \dot{e} &= g(x, e)\end{aligned}\right\} (x, e) \in C, \quad \left.\begin{aligned}x^+ &= x \\ e^+ &= 0\end{aligned}\right\} (x, e) \in D, \quad (6)$$

---
[1] A set $\mathcal{A} \subset \mathbb{R}^n$ is forward invariant for system (1), if for any solution $\phi$ to (1) with $\phi(0, 0) \in \mathcal{A}$ we have $\phi(t, j) \in \mathcal{A}$ for all $(t, j) \in \text{dom}\,\phi$.

where $x = (x_P, x_C) \in \mathbb{R}^{n_x}$, $e = (e_{x_P}, e_u) \in \mathbb{R}^{n_e}$ denotes the sampling-induced error with $e_{x_P} = \hat{x}_P - x_P$, $e_u = \hat{u} - u$. The sets $C$ and $D$ are closed and respectively denote the flow and the jump sets, they are defined according to the triggering condition. Typically, the system flows on $C$ and experiences a jump on $D$ where the triggering condition is satisfied. When $(x,e) \in C \cup D$, the system can either jump or flow, the latter only if flowing keeps $(x,e)$ in $C$. Functions $f$ and $g$ are defined as (where we can replace $\hat{x}_P$ by $x_P + e_{x_P}$):

$$f(x,e) = \begin{pmatrix} f_P(x_P, g_C(x_C, \hat{x}_P) + e_u) \\ f_C(x_C, \hat{x}_P) \end{pmatrix}$$

$$g(x,e) = \begin{pmatrix} \hat{f}_P(x_P, x_C, \hat{x}_P, g_C(x_C, \hat{x}_P) + e_u) \\ -f_P(x_P, g_C(x_C, \hat{x}_P) + e_u) \\ \hat{f}_C(x_P, x_C, \hat{x}_P, g_C(x_C, \hat{x}_P) + e_u) \\ -\frac{\partial g_C}{\partial x_C}(x_C, \hat{x}_P) f_C(x_C, \hat{x}_P) \\ -\frac{\partial g_C}{\partial \hat{x}_P}(x_C, \hat{x}_P) \\ \times \hat{f}_P(x_P, x_C, \hat{x}_P, g_C(x_C, \hat{x}_P) + e_u) \end{pmatrix}$$

(7)

and are assumed to be continuous.

**Remark 2.** *Our assumptions allow for triggering rules that depend both on $x$ and $e$. However, the specific choice of triggering rule needs to be done according to the implementation scenario. In the case of dynamic controllers, a triggering rule depending on $x_C$ requires continuous communication between the sensors and the controller. This is difficult to achieve in practice since sensors do not have, in general, access to the state of the controller. We have chosen to present the problem in a general setting because it allows to recover as particular cases the stabilization using a static controller (as in Sections III-B, V-A for example) and the cases where only the plant states or the inputs are sampled.*

The main problem addressed in this paper is to define the triggering condition, i.e. the flow and jump sets $C$ and $D$ in (6), in order to minimize the resource usage while ensuring asymptotic stability properties. We now introduce the main idea of the framework presented hereafter in Section IV by interpreting the work in [19] using a hybrid Lyapunov function.

*B. Main idea*

We first revisit the work in [19] where a static controller $u = g_C(x_P)$ is assumed to render the closed-loop system (3) input-to-state stable (ISS) with respect to the sampling-induced errors (that can be considered as measurement errors at this stage since, when the controller is static, the sampling-induced error can be seen to be only due to the sampling of the measurements, i.e. $e = e_{x_P}$). This is equivalent to the following assumption (see Theorem 1 in [18]) where $x = x_P$ (as the controller is static).

**Assumption 1.** *There exists a smooth Lyapunov function $V : \mathbb{R}^{n_x} \to \mathbb{R}$ and $\underline{\alpha}_V, \overline{\alpha}_V, \alpha, \gamma \in \mathcal{K}_\infty$ such that for all $x \in \mathbb{R}^{n_x}$:*

$$\underline{\alpha}_V(|x|) \leq V(x) \leq \overline{\alpha}_V(|x|), \quad (8)$$

*and for all $(x,e) \in \mathbb{R}^{n_x + n_e}$:*

$$\frac{\partial V}{\partial x} f(x,e) \leq -\alpha(V(x)) + \gamma(|e|). \quad (9)$$

Since zero-order-hold devices are used in [19], we have $g(x,e) = -f(x,e)$ in (7) and the model (6) is here:

$$\left. \begin{array}{l} \dot{x} = f(x,e) \\ \dot{e} = -f(x,e) \end{array} \right\} (x,e) \in C, \quad \left. \begin{array}{l} x^+ = x \\ e^+ = 0 \end{array} \right\} (x,e) \in D. \quad (10)$$

From (9), we deduce that $\sigma \alpha(V(x)) \geq \gamma(|e|)$ with $\sigma \in (0,1)$ implies:

$$\frac{\partial V}{\partial x} f(x,e) \leq -(1-\sigma)\alpha(V(x)). \quad (11)$$

In that way, the triggering rule in [19] can be written as $\sigma \alpha(V(x)) \leq \gamma(|e|)$, that we rewrite as:

$$V(x) \leq \alpha^{-1}(\sigma^{-1}\gamma(|e|)) =: \tilde{\gamma}(|e|). \quad (12)$$

At each transmission instant, $e$ is reset to 0, so we have $\tilde{\gamma}(|e^+|) = 0 \leq V(x)$ and $V$ decreases monotonically according to (11). The next transmission occurs as soon as (12) is satisfied. The flow and the jump sets in (10) can be defined as follows:

$$\begin{array}{l} C = \{(x,e) : \tilde{\gamma}(|e|) \leq V(x)\} \\ D = \{(x,e) : \tilde{\gamma}(|e|) \geq V(x)\}. \end{array} \quad (13)$$

To guarantee the existence of a minimum interval of time between two transmissions when $(x,e) \neq 0$, the following conditions are used in [19].

**Assumption 2.** *For any compact set $S \subset \mathbb{R}^{n_x + n_e}$, there exist $L_1, L_2 \in \mathbb{R}_{\geq 0}$ such that for all $(x,e) \in S$:*

$$|f(x,e)| \leq L_1(|x| + |e|) \quad (14)$$

$$\underline{\alpha}_V^{-1} \circ \tilde{\gamma}(|e|) \leq L_2|e|. \quad (15)$$

The stability analysis of system (10) can be done using the following Lyapunov function (assuming $\tilde{\gamma}$ is locally Lipschitz): $R(x,e) = \max\{V(x), \tilde{\gamma}(|e|)\}$. Indeed,

- *Property (a):* $R$ is positive definite and radially unbounded in view of (8) and since $\tilde{\gamma} \in \mathcal{K}_\infty$.
- *Property (b):* $R$ decreases on $C$ according to (11).
- *Property (c):* $R$ does not increase at jumps since $x^+ = x$ and $e^+ = 0$.
- *Property (d):* it was shown in [19] using Assumption 2 that there does exist a uniform minimal time interval between two successive transmission instants (as long as $(x,e) \neq 0$) for solutions that start in a compact set that contains the origin. In other words, solutions to (10) have a semiglobal dwell time on $\mathbb{R}^{n_x+n_e} \setminus \{0\}$ according to Definition 2.

We show in Section IV that these four properties guarantee asymptotic stability properties for system (10) and that they can be used to build up other event-triggering conditions. Note that similar ingredients are used to prove the stability of other types of hybrid systems in [13], [14] for example.

IV. A LYAPUNOV-BASED FRAMEWORK

Before stating the main result of this section, it is important to note that auxiliary variables may be introduced to define the triggering condition. Indeed, it is common in the hybrid literature to introduce additional variables like clocks to ensure or analyse the stability e.g. [5], [13]. We will see in Section V-A that the strategy in [20] can be interpreted using our framework by making use of an additional variable which is employed to build up a decreasing threshold on the known Lyapunov function for the system in the absence of sampling. We also show in Section V-B that the event-triggered policy in Section III-B can be redesigned to exhibit larger inter-event intervals thanks to the use of an auxiliary variable. Therefore, we denote by a single vector variable $\eta \in \mathbb{R}^{n_\eta}$ the additional variables which may be needed for describing the system that are neither $x$ nor $e$.

In that way, to define a triggering condition ensuring desired stability properties for the overall system is tantamount to defining appropriate flow and jump sets $C$ and $D$ for the following hybrid system:

$$\left.\begin{array}{rcl}\dot{x} & = & f(x,e) \\ \dot{e} & = & g(x,e) \\ \dot{\eta} & = & h(x,e,\eta)\end{array}\right\} q \in C, \quad \left.\begin{array}{rcl}x^+ & = & x \\ e^+ & = & 0 \\ \eta^+ & = & l(x,e,\eta)\end{array}\right\} q \in D, \quad (16)$$

where $q = (x, e, \eta) \in \mathbb{R}^{n_q}$, $h, l$ are continuous and $C, D$ are closed subsets of $\mathbb{R}^{n_q}$. We use $\dot{q} = F(q)$ and $q^+ = G(q)$ to denote (16).

The stability of system (16) can be guaranteed by means of the following theorem. It can be seen as a variation of the results in [5].

**Theorem 1.** *Consider system (16) and suppose $G(D) \subset (C \cup D)$ and that there exist a locally Lipschitz function $R : \mathbb{R}^{n_q} \to \mathbb{R}$ and a continuous function $\upsilon : \mathbb{R}^{n_\eta} \to \mathbb{R}^{n_\upsilon}$ with $n_\upsilon \leq n_\eta$ such that the following conditions hold.*

(i) *There exist $\underline{\alpha}_R, \overline{\alpha}_R \in \mathcal{K}_\infty$ such that for any $q \in C \cup D$, $\underline{\alpha}_R(|(x, e, \upsilon(\eta))|) \leq R(q) \leq \overline{\alpha}_R(|(x, e, \upsilon(\eta))|)$.*
(ii) *There exists $\alpha_R \in \mathcal{K}_\infty$ such that for all $q \in C$: $R^\circ(q; F(q)) \leq -\alpha_R(R(q))$.*
(iii) *For all $q \in D$, $R(G(q)) \leq R(q)$.*
(iv) *Solutions to (16) have a semiglobal dwell time on $\mathbb{R}^{n_q} \setminus \mathcal{A}$, where $\mathcal{A} = \{q : (x, e, \upsilon(\eta)) = 0\}$.*

*Then the set $\mathcal{A}$ is S-GAS.*

Theorem 1 provides a Lyapunov-based prescriptive framework for developing event-triggered control strategies for nonlinear systems as we show in Section V. Other triggering rules may be derived by following the guidelines below for instance. We illustrate each item with the example of Section III-B for the sake of clarity.

1) Select a locally Lipschitz function $R : \mathbb{R}^{n_q} \to \mathbb{R}$ that satisfies item (i) of Theorem 1. Usually, $R$ is built using a known Lyapunov function $V$ for the continuous-time system (3)-(4) in the absence of sampling and a positive definite radially unbounded function $W(e)$ that has to be designed. Typically, $W$ is chosen by investigating the robustness property of the closed-loop system $\dot{x} = f(x, e)$ with respect to $e$ that is assumed to hold. The sets $C$ and $D$ have not been defined so far but item (i) of Theorem 1 needs to hold on $C \cup D$. This apparent contradiction is overcome as follows. When there is no variable $\eta$, as it is the case so far, we typically have $C \cup D = \mathbb{R}^{n_x + n_e}$ and we do not need to know $C$ and $D$ to verify item (i) of Theorem 1.

   Section III-B: we took $W(e) = \tilde{\gamma}(|e|)$ where $\tilde{\gamma}(|e|)$ is defined in (12) that is deduced from the ISS property stated in Assumption 1. We considered $R(q) = \max\{V(x), \tilde{\gamma}(|e|)\}$ with $q = (x, e)$, that does satisfy item (i) of Theorem 1 on $\mathbb{R}^{n_x + n_e}$. This corresponds to Property (a).

2) Choose $\alpha_R \in \mathcal{K}_\infty$ for item (ii) of Theorem 1. Obviously, if $R(x, 0) = V(x)$, the decreasing rate $\alpha_R$ will have to be less than the decreasing rate of $V$ in order to allow some flow before entering the set $D$.
   Section III-B: we have taken $\alpha_R(s) = (1 - \sigma)\alpha(s) \leq \alpha(s)$ for $s \in \mathbb{R}_{\geq 0}$, since $\sigma \in (0, 1)$.

3) Define the flow and the jump sets to be closed and such that items (i)-(iii) of Theorem 1 hold and $G(D) \subset (C \cup D)$. For instance, when items (i) and (iii) of Theorem 1 are satisfied for all $q \in \mathbb{R}^{n_q}$, we can directly take the following sets: $C = \{q : R^\circ(q; F(q)) \leq -\alpha_R(R(q))\}$ and $D = \{q : R^\circ(q; F(q)) \geq -\alpha_R(R(q))\}$ which ensure item (ii) of Theorem 1 and $G(D) \subset (C \cup D) = \mathbb{R}^{n_q}$.
   Section III-B: the flow and the jump sets in (13) guarantee that items (ii)-(iii) of Theorem 1 holds in view of (11) and since $x^+ = x$ and $e^+ = 0$, that is equivalent to Properties (b)-(c). We note that $C \cup D = \mathbb{R}^{n_x + n_e}$.

4) Study the existence of dwell times. Among other techniques, Lemma 1 below can be used for this purpose. The existence of dwell times notably depends on the triggering condition and the vector field $F$ that is usually assumed to satisfy some Lipschitz properties. If the existence of a dwell time is guaranteed, the desired result is obtained. Otherwise, variable $\eta$ may be introduced, then go back to 1) and modify the function[2] $R$. The way the variable $\eta$ may be chosen will become clearer in the light of Section V. The non-existence of dwell time may also be due to the fact that the decreasing rate of $R$ along flows, $\alpha_R$ in 2), is too strong, thus choose a different function $\tilde{\alpha}_R \in \mathcal{K}_\infty$ such that $\tilde{\alpha}_R(s) < \alpha_R(s)$ for any $s \in \mathbb{R}_{>0}$.
   Section III-B: the existence of semiglobal dwell-time solutions is guaranteed in [19] using Assumption 2, as stated in Property (d).

The following lemma provides a tool for verifying the existence of dwell times which is used in the proofs of the theorems of Section V.

**Lemma 1.** *Consider system (16) and suppose the following holds.*

(i) *$G(D) \subset (C \cup D)$ and items (i)-(iii) of Theorem 1 are satisfied.*
(ii) *For any $q \in D$, $G(0, e, \eta) \in \mathcal{A} \cup \mathcal{M}$ where $\mathcal{A} = \{q : (x, e, \upsilon(\eta)) = 0\}$ and $\mathcal{M} \subset (C \setminus D)$ is a forward invariant set.*
(iii) *There exists a locally Lipschitz function $\psi : \Theta(\mu) \to$*

---
[2]It may be the case that the new $R$ no longer satisfies item (i) of Theorem 1 on $\mathbb{R}^{n_q}$, so identify $S \subseteq \mathbb{R}^{n_\eta}$ such that item (i) of Theorem 1 is satisfied on $\mathbb{R}^{n_x + n_e} \times S$, afterwards make sure $(C \cup D) \subset \mathbb{R}^{n_x + n_e} \times S$.

$\mathbb{R}_{\geq 0}$, where $\Theta(\mu) = \{q \in C \cup D : R(q) \leq \mu \text{ and } x \neq 0\}$ for $\mu > 0$, such that:

(iii-a) There exists $a \in \mathbb{R}_{\geq 0}$ such that for any $q \in D$ with $G(q) \in \Theta(\mu)$: $\psi(G(q)) \leq a$.

(iii-b) There exists $b > a$ such that for any solution $\phi$ to (16), $t_j \leq t$ with $(t,j) \in \text{dom}\,\phi$: $(\psi(\phi(t,j)) < b) \Rightarrow (\phi(t,j) \in C \backslash D)$.

(iii-c) There exists a continuous non-decreasing function $\lambda : \mathbb{R}_{>0} \to \mathbb{R}_{\geq 0}$ such that for all $q \in \Theta(\mu)$: $\psi^\circ(q; F(q)) \leq \lambda(\psi(q))$.

Then solutions to (16) have a semiglobal dwell time on $\mathbb{R}^{n_q} \backslash \mathcal{A}$.

The conditions of Lemma 1 can be interpreted as follows. Item (i) simply states that all the conditions of Theorem 1 are verified except item (iv) which is the purpose of this lemma. When a jump occurs when $x = 0$, we know from (16) that after the jump $(x, e) = 0$, then we no longer need to transmit since 0 is very likely to be an equilibrium point of the $(x, e)$-system. That is what item (ii) of Lemma 1 says: after a jump when $x = 0$, solutions to (25) lies in the stable set $\mathcal{A}$ or in a subset of $C \backslash D$ and will never leave it. The function $\psi$ which is considered in item (iii) of Lemma 1 is used to guarantee that there exists a minimum uniform amount of time between two successive jumps. By estimating the time it takes for $\psi$ to grow from $a$ to $b$, we are able to obtain the desired result.

**Remark 3.** *The triggering condition that satisfies the conditions of Theorem 1 respects the practical requirement that there does exist a uniform minimum time interval between two transmissions according to item (iv) of Theorem 1. The only region of the state space where this may not be guaranteed is when $(x, e, \upsilon(\eta)) = 0$, but this will only occur if the system is initialized in the stable set.*

## V. APPLICATIONS

We already know that the framework allows us to capture the work in [19], we show in this section that it is also the case for the strategy in [20]. Afterwards, new triggering rules are proposed.

### A. Event-triggered strategy in [20]

As in Section III-B, the controller is static ($x = x_P$) and implemented using zero-order-hold devices. It is considered that Assumption 1 is satisfied with $\alpha$ linear that is $\alpha(s) = \bar{\alpha}s$ with $\bar{\alpha} \in \mathbb{R}_{>0}$. The triggering rule is defined to guarantee that $V(x(t_j))$ always decreases at a certain rate compared to $V(x(t_{j-1}))$. In that way, the control loop is closed in [20] as soon as the condition below is violated, for $t \geq t_{j-1}$:

$$V(x(t)) \leq -\bar{\sigma}\bar{\alpha}V(x(t_{j-1}))(t - t_{j-1}) + V(x(t_{j-1})) \\ = (-\bar{\sigma}\bar{\alpha}(t - t_{j-1}) + 1)V(x(t_{j-1})), \quad (17)$$

where $\bar{\sigma} \in (0, 1)$. Since zero-order-hold devices are considered, we have $\hat{x}(t) = x(t_{j-1})$ and $e(t) = e_{x_P}(t) = x(t_{j-1}) - x(t)$. Consequently, (17) is equivalent to, for $t \geq t_{j-1}$:

$$V(x(t)) \leq (-\bar{\sigma}\bar{\alpha}(t - t_{j-1}) + 1)V(\hat{x}(t) + x(t) - x(t)) \\ = (-\bar{\sigma}\bar{\alpha}(t - t_{j-1}) + 1)V(x(t) + e(t)). \quad (18)$$

To model (17) using the hybrid formulation (16), we introduce the variable $\eta \in \mathbb{R}$ as the solution of $\dot{\eta} = -\bar{\sigma}\bar{\alpha}$ on flows and $\eta^+ = 1$ at jumps. We see that $\eta(t) = -\bar{\sigma}\bar{\alpha}(t - t_{j-1}) + 1$ for $t \in [t_{j-1}, t_j]$ ($j \in \mathbb{Z}_{>0}$). In that way, we can reformulate (18) using the following algebraic inequality:

$$V(x) \leq \eta V(x + e). \quad (19)$$

The problem can then be modeled as follows:

$$\left.\begin{array}{l}\dot{x} = f(x, e) \\ \dot{e} = -f(x, e) \\ \dot{\eta} = -\bar{\sigma}\bar{\alpha}\end{array}\right\} q \in C, \quad \left.\begin{array}{l}x^+ = x \\ e^+ = 0 \\ \eta^+ = 1\end{array}\right\} q \in D, \quad (20)$$

where $q = (x, e, \eta)$,

$$\begin{array}{l}C = \{q : V(x) \leq \eta V(x + e) \text{ and } \eta \in [\varepsilon, 1]\} \\ D = D_1 \cup D_2,\end{array} \quad (21)$$

with $D_1 = \{q : V(x) \geq \eta V(x+e) \text{ and } \frac{\partial V}{\partial x}f(x,e) \geq -\bar{\sigma}\bar{\alpha}V(x)\}$ and $D_2 = \{q : \eta = \varepsilon\}$ where $\varepsilon \in (0, 1)$ is arbitrary small. The condition $\frac{\partial V}{\partial x}f(x,e) \geq -\bar{\sigma}\bar{\alpha}V(x)$ has been added in the definition of $D_1$ to avoid Zeno behaviour since after a jump $V(x) = \eta V(x + e)$ holds. Indeed, it is not necessary to jump again since $V(x)$ will decrease faster than $\eta V(x + e)$ for some time according to (9). The lower bound $\varepsilon$ on $\eta$ is used to guarantee that the threshold on $V(x)$ defined by $\eta V(x+e)$ (see (19)) never reaches the origin when $V(x+e) \neq 0$. This condition adds no conservatism as by setting $\varepsilon$ sufficiently small, the triggering condition $\eta = \varepsilon$ will not be satisfied in practice before $q$ reaches $D_1$.

We recover Theorem 3.2 in [20] and relax some of the required conditions.

**Theorem 2.** *Consider system (20) and suppose Assumption 1 holds with $\alpha(s) = \bar{\alpha}s$ ($\bar{\alpha} \in \mathbb{R}_{>0}$) and Assumption 2 is satisfied. Then the set $\mathcal{A} = \{q : (x, e) = 0\}$ is S-GAS and solutions to (20) have a semiglobal dwell time on $\mathbb{R}^{n_q} \backslash \mathcal{A}$.*

We note that the conditions of Theorem 2 are more general than those of Theorem 3.2 in [20] as $\gamma$ in (9) is allowed to be a nonlinear function. In addition, condition (15) in this paper extends (5) in [20] and allows us to consider more general types of Lyapunov functions, such as quadratic, which is not the case in [20].

### B. New triggering rules

In Section V-A, the triggering condition is obtained by defining a decreasing threshold on $V$ (see (17)). In this subsection, we propose an alternative that consists in defining a similar threshold for an appropriate function $W$ for the $e$−system. We suppose that the dynamic controller (4) has been designed so that Assumption 1 applies. Thus, by using the ISS property of the $x$-system, we will be able to show that when $W$ remains below a given decreasing threshold, system (16) satisfies asymptotic stability properties.

We define our threshold variable $\eta \in \mathbb{R}$ as the solution of the following differential equation on flows:

$$\dot{\eta} = -\delta(\eta), \quad (22)$$

where $\delta$ is any class-$\mathcal{K}_\infty$ function, and at jumps,

$$\eta^+ = \tilde{\gamma}(|e|) =: W(e), \tag{23}$$

where $\tilde{\gamma}(s) = \alpha^{-1}(\sigma^{-1}\gamma(s))$ for $s \in \mathbb{R}_{\geq 0}$, with $\sigma \in (0,1)$ as in (12). We note that $W$ is positive definite and radially unbounded. An obvious choice of triggering rule is: $W(e) \geq \eta$. Nevertheless, in the case where $W(e) \leq V(x)$, $V$ decreases according to (11) and therefore we do not need to close the loop. This suggests considering the following triggering condition instead:

$$W(e) \geq \max\{\eta, V(x)\}. \tag{24}$$

The problem can be modeled as follows:

$$\left.\begin{array}{l}\dot{x} = f(x,e) \\ \dot{e} = g(x,e) \\ \dot{\eta} = -\delta(\eta)\end{array}\right\} q \in C, \quad \left.\begin{array}{l}x^+ = x \\ e^+ = 0 \\ \eta^+ = W(e)\end{array}\right\} q \in D, \tag{25}$$

where $q = (x, e, \eta)$,

$$\begin{array}{l}C = \{q : \max\{V(x), \eta\} \geq W(e) \text{ and } \eta \geq 0\} \\ D = \{q : \max\{V(x), \eta\} \leq W(e) \text{ and } \eta \geq 0\}.\end{array} \tag{26}$$

The following theorem ensures the stability of system (25).

**Theorem 3.** *Consider system (25), suppose the following conditions hold.*
  (i) *Assumptions 1-2 apply.*
  (ii) *Function $\tilde{\gamma}$ is locally Lipschitz.*
  (iii) *For any compact set $S \subset \mathbb{R}^{n_x+n_e}$, there exist $L_3 \in \mathbb{R}_{\geq 0}$ such that for all $(x, e) \in S$: $|g(x, e)| \leq L_3(|x| + |e|)$.*

*Then $q = (x, e, \eta) = 0$ is S-GAS and solutions to (25) have a semiglobal dwell time on $\mathbb{R}^{n_q}\setminus\{0\}$.*

Contrary to Section V-A, we note that $\alpha$ in (9) is allowed to be nonlinear. In addition, we do not focus on zero-order-hold devices that is why condition (iii) of Theorem 3 is introduced in order to guarantee the existence of dwell times. We show on an example in Section VI that the inter-event intervals can be enlarged to some extent compared to Section V-A by playing with the initial value of $\eta$.

## VI. ILLUSTRATIVE EXAMPLE

To illustrate the benefits of the strategy presented in Section V-B, we revisit the example considered in [13]. The simplified version of the considered nonlinear system is:

$$\dot{x} = dx^2 - x^3 + u, \tag{27}$$

where $d$ is an unknown possibly time-varying parameter satisfying $|d| < 1$. The stabilizing control law considered in [13] was $u = -2x$. We select $V(x) = \frac{1}{2}x^2$ as our Lyapunov function that satisfies Assumption 1 with[3] $\underline{\alpha}_V(s) = \overline{\alpha}_V(s) = \frac{1}{2}s^2$, $\alpha(s) = 0.84s$ and $\gamma(s) = 2.66s^2$ for $s \in \mathbb{R}_{\geq 0}$. We consider 200 random initial conditions distributed in the interval $[-1, 1]$. The parameter $d$ takes for each initial condition a random value in the interval $[0, 1]$. We compare the average number of executions required under the technique in [20] as extended in Section V-A, the event-triggered strategy proposed in Section V-B and the periodic strategy in [13] in Table I, for different values for the design parameter $\eta(0, 0)$. We select $\bar{\sigma} = 10^{-3}$ in (20) and $\delta(s) = 0.5s$ with $s \in \mathbb{R}_{\geq 0}$ in (22). It can be observed that the average number of executions is considerably lower under the event-triggered strategies. Moreover we note that the policy in Section V-B generates less executions than [20] and that it can be adjusted by means of the design parameter $\eta(0, 0)$. It is however not easy to compare the performance guarantees under the two different policies, since Theorem 2 and Theorem 3 ensure different stability properties.

[3]The Yalmip software ([12]) was used to compute $\alpha$ and $\gamma$.

| [13] | [20] | Section V-B | | |
| | | $\eta(0,0) = 0.1$ | $\eta(0,0) = 1$ | $\eta(0,0) = 2$ |
| --- | --- | --- | --- | --- |
| 54.34 | 18.47 | 16.88 | 15.27 | 13.31 |

TABLE I
AVERAGE NUMBER OF EXECUTIONS OVER 200 INITIAL CONDITIONS
FOR A SIMULATION TIME OF 20S FOR THE EXAMPLE IN [13].

We can explore deeper the role played by the initial condition of the auxiliary variable $\eta(0, 0)$ and its effect on performance and number of executions. Simulations have shown that smaller values of $\eta(0, 0)$ imply a faster decay, at the expense of more executions. Hence the design parameter $\eta(0, 0)$ represents the tradeoff between performance and resource usage. Similar conclusions can be drawn for the decay function $\delta$ in (22).

The technique in Section V-B exhibits great potential for real-time scheduling, since both the initial value for the auxiliary variable and the differential equation in (22) can be designed according to the available resources. For instance, functions $\delta$ with slow increasing slopes could be chosen in case of overload in the network or in the processor executing the controller.

## VII. APPENDIX

**Proof of Theorem 1.** Let $B = \{q \in C \cup D : |q|_\mathcal{A} < \Delta\}$ where $\Delta \in \mathbb{R}_{>0}$ and $\phi$ be a solution to (16) with $\phi(0, 0) \in B$. Define the set $\Omega = \{q \in C \cup D : R(q) \leq \mu\}$ where $\mu \in \mathbb{R}_{>0}$ is such that $B \subseteq \Omega$ (take for instance $\mu = \overline{\alpha}_R(\Delta)$ in view of item (i) of Theorem 1). The set $C \cup D$ is forward invariant for system (16) since $G(D) \subset (C \cup D)$. Hence, in view of items (i)-(iii) of Theorem 1, $\phi(t, j) \in \Omega$ for any $(t, j) \in \text{dom}\,\phi$. From item (ii) of Theorem 1 and by using standard comparison principles, there exists $\beta \in \mathcal{KL}$ that satisfies, for all $(s, t_1, t_2) \in \mathbb{R}_{\geq 0} \times \mathbb{R}_{\geq 0} \times \mathbb{R}_{\geq 0}$:

$$\beta(s, t_1 + t_2) = \beta(\beta(s, t_1), t_2), \tag{28}$$

and such that, for all $(t_j, j) \preceq (t, j) \in \text{dom}\,\phi$,

$$R(\phi(t, j)) \leq \beta(R(\phi(t_j, j)), t - t_j), \tag{29}$$

where $(t_j, j) \preceq (t, j)$ means that $t_j \leq t$. From item (iii) of Theorem 1, it follows that:

$$R(\phi(t_{j+1}, j+1)) \leq R(\phi(t_{j+1}, j)) \tag{30}$$

for all $j$ such that $(t, j) \in \text{dom}\, \phi$ for some $t \in \mathbb{R}_{\geq 0}$. Combining (28)-(30), we obtain:

$$R(\phi(t,j)) \leq \beta(R(\phi(0,0)), t) \qquad \forall (t,j) \in \text{dom}\, \phi. \quad (31)$$

Now let $(t,j) \in \text{dom}\, \phi$, if $\phi(t,j) \in \mathcal{A}$ then $R(\phi(t,j)) = 0$ according to item (i) of Theorem 1. If $\phi(t,j) \notin \mathcal{A}$ then that means that $\phi(t', j') \notin \mathcal{A}$ for all $(t', j') \preceq (t,j) \in \text{dom}\, \phi$ (i.e. $t' \leq t$ and $j' \leq j$) since $\mathcal{A}$ is forward invariant for system (16) (according to items (i)-(iii) of Theorem 1 and since $G(D) \subset (C \cup D)$). As a consequence, we have that $t' \geq \tau(\mu) j'$ for all $(t', j') \preceq (t,j) \in \text{dom}\, \phi$, where $\tau(\Delta) \in \mathbb{R}_{>0}$ is a minimal interval of times between two jumps on $\Omega \backslash \mathcal{A}$ whose existence is ensured by item (iv) of Theorem 1. It follows that:

$$R(\phi(t,j)) \leq \beta(R(\phi(0,0)), \tfrac{1}{2}t + \tfrac{1}{2}\tau(\Delta)j). \quad (32)$$

By using item (i) of Theorem 1, we deduce that, for all $(t,j) \in \text{dom}\, \phi$:

$$|\phi(t,j)|_{\mathcal{A}} \leq \underline{\alpha}_R^{-1}\Big(\beta\big(\overline{\alpha}_R(|\phi(0,0)|_{\mathcal{A}}), \tfrac{1}{2}t + \tfrac{1}{2}\tau(\Delta)j\big)\Big), \quad (33)$$

denoting $\beta_\Delta : (s,t,j) \mapsto \underline{\alpha}_R^{-1}\Big(\beta\big(\overline{\alpha}_R(s), \tfrac{1}{2}t + \tfrac{1}{2}\tau(\Delta)j\big)\Big) \in \mathcal{KLL}$ (since $\beta \in \mathcal{KL}$, $\underline{\alpha}_R, \overline{\alpha}_R \in \mathcal{K}_\infty$), for all $(t,j) \in \text{dom}\, \phi$:

$$|\phi(t,j)|_{\mathcal{A}} \leq \beta_\Delta(|\phi(0,0)|_{\mathcal{A}}, t, j). \quad (34)$$

Hence, the set $\mathcal{A}$ is S-GAS according to Definition 1. $\square$

**Proof of Lemma 1.** Let $\Delta \in \mathbb{R}_{>0}$ and define $B = \{q \in C \cup D : |q|_{\mathcal{A}} < \Delta\}$, $\Omega = \{q \in C \cup D : R(q) \leq \mu\}$ where $\mu \in \mathbb{R}_{>0}$ is such that $B \subseteq \Omega$ as in the proof of Theorem 1. Let $\phi$ be a solution to (16) with $\phi(0,0) \in B$. According to items (i)-(iii) of Theorem 1, $\phi(t,j) \in \Omega$ for any $(t,j) \in \text{dom}\, \phi$. Denote $t_1 \in (0, \infty)$ the first jump instant with $(t_1, 1) \in \text{dom}\, \phi$ (if no jump ever occurs i.e. $t_1 = \infty$, (2) is obviously satisfied). We have that $\phi(t_1, 0) \in D$ and $R(\phi(t_1, 1)) \leq \mu$ in view of item (iii) of Theorem 1. If $\phi(t_1, 1) = G(\phi(t_1, 0)) \notin \Theta(\mu)$ then $x(t_1, 0) = x(t_1, 1) = 0$ and $\phi(t_1, 1) \in \mathcal{A} \cup \mathcal{M}$ since it is assumed that $G(0, e, \eta) \in \mathcal{A} \cup \mathcal{M}$ for any $q \in D$. As a consequence, $\phi(t, j) \in \mathcal{A} \cup \mathcal{M}$ for any $(t_1, 1) \preceq (t,j) \in \text{dom}\, \phi$ as $\mathcal{A}$ and $\mathcal{M}$ are forward invariant for system (16). In that way, (2) is ensured since on $\mathcal{M} \subset C \backslash D$ no jump ever occurs. If $\phi(t_1, 1) \in \Theta(\mu)$, then according to item (iii-a) of Lemma 1, $\psi(G(\phi(t_1, 0))) = \psi(\phi(t_1, 1)) \leq a < b$ and therefore a jump cannot occur immediately in view of item (iii-b) of Lemma 1. Let denote $t_2 > t_1$ with $(t_2, 2) \in \text{dom}\, \phi$ the next jump instant and suppose $t_2 \neq \infty$, otherwise the desired result holds. By the continuity of $\psi$ and the solution $\phi$ to system (16) on flows, there exists $t^* \in (t_1, t_2]$ such that $\psi(\phi(t^*, 1)) = b$. According to item (iii-b) of Lemma 1, we deduce that $\phi(t, 1) \in C \backslash D$ for any $t \in [t_1, t^*)$. In view of item (iii-c) of Lemma 1, invoking standard comparison principles, we deduce that $\psi(\phi(t, 1)) \leq \theta(t)$ for any $t \in [t_1, t_2)$ where $\theta(t)$ is the solution of $\dot\theta = \lambda(\theta)$ satisfying $\theta(t_1) = a \geq \psi(\phi(t_1, 1))$. The next jump cannot occur before the time $\tau(\Delta) \in (0, \infty]$ it takes for $\theta$ to evolve from $a$ to $b$ (which is independent of $t_1$) has elapsed. Note that $\psi(\phi(\cdot, 1))$ cannot reach $b$ before $\theta$. By induction, we deduce that the inter-jump interval on $B \backslash \mathcal{A}$ is lower bounded by $\tau(\Delta)$. Hence solutions to (16) have a semiglobal dwell time on $\mathbb{R}^{n_q} \backslash \mathcal{A}$ according to Definition 2.$\square$

**Proof of Theorem 2.** The proof consists in applying Theorem 1 for system (20). Let $q = (x, e, \eta) \in D$, we know that $G(q) = (x, 0, 1)$ therefore $G(q) \in C$ and we deduce that $G(D) \subset C \cup D$. Consider system (20) and the following candidate Lyapunov function for $q \in \mathbb{R}^{n_q}$:

$$R(q) = \max\{V(x), \eta V(x+e)\}. \quad (35)$$

We first prove that $R$ satisfies item (i) of Theorem 1 with $\nu(\eta) = 0$. Let $q \in C \cup D$ then $\eta \in [\varepsilon, 1]$. We see that $R(0, 0, \eta) = 0$. Now suppose $R(q) = 0$, since $R(q) \geq V(x) \geq 0$, we have $V(x) = 0$, so $x = 0$ according to (8). Similarly $R(q) \geq \eta V(x+e) \geq \varepsilon V(x+e) \geq 0$ (since $\eta \geq \varepsilon$), thus $x + e = 0$ but since $x = 0$ we obtain $e = 0$. Hence $R(q) = 0$ if and only if $(x, e) = 0$. We now verify that $R(q) \to \infty$ when $|(x,e)| \to \infty$. Noting that $R(q) \geq V(x)$, we see that $R(q) \to \infty$ as $|x| \to \infty$ from (8). Suppose $|e| \to \infty$ and $|x|$ does not tend to $\infty$, without loss of generality we fix $x$. Take any sequence $\{e_i\}_{i \in \mathbb{Z}_{>0}}$ such that $|e_i| \to \infty$ as $i \to \infty$. In view of (35), there exists $i_0$ such that $R(q) = \eta V(x + e_i)$ for any $i \geq i_0$ since $|x|$ does not tend to $\infty$. Therefore, $R(q) \to \infty$ as $i \to \infty$ in view of (8). Since $\{e_i\}_{i \in \mathbb{Z}_{>0}}$ has been chosen arbitrarily, $R(q) \to \infty$ as $|e| \to \infty$. We have shown that $R(q) \to \infty$ when $|(x,e)| \to \infty$. Consequently, there exist $\underline{\alpha}_R, \overline{\alpha}_R \in \mathcal{K}_\infty$ so that item (i) of Theorem 1 applies by following similar lines as in the proof of Lemma 4.3 in [9] and using the fact that $\eta \in [\varepsilon, 1]$. On $C$, we have that $V(x) \leq \eta V(x+e)$ so $R(q) = \eta V(x+e)$. Then, $R^\circ(q; F(q)) = -\bar\sigma \bar\alpha V(x+e) + \eta \frac{\partial V}{\partial x}(f(x+e) - f(x+e)) = -\bar\sigma \bar\alpha V(x+e)$. Since $\eta \leq 1$, $R^\circ(q; F(q)) \leq -\bar\sigma \bar\alpha \eta V(x+e) = -\bar\sigma \bar\alpha R(q)$ and item (ii) of Theorem 1 holds with $\alpha_R(s) = \bar\sigma \bar\alpha s$ for $s \in \mathbb{R}_{\geq 0}$. Let $q \in D$, $R(G(q)) = \max\{V(x), V(x+0)\} \leq R(q)$, item (iii) of Theorem 1 is ensured. We now show that item (iv) of Theorem 1 is satisfied by applying Lemma 1. First, we note that $G(0, e, \eta) = (0, 0, 1) \in \mathcal{A}$ for $q \in D$, where $\mathcal{A} = \{q : (x, e) = 0\}$, so that item (ii) of Lemma 1 holds. We take:

$$\psi : q \mapsto \max\{L_2 \tfrac{|e|}{|x|}, \tfrac{1-\eta}{1-\varepsilon}\} \quad (36)$$

that is defined on $\Theta(\mu)$ (see Lemma 1) with $\mu \in \mathbb{R}_{>0}$ and $L_2$ comes from Assumption 2. We see that for any $q \in D$ with $G(q) \in \Theta(\mu)$, $\psi(G(q)) = 0$: item (iii-a) of Lemma 1 is ensured with $a = 0$. Let $q \in \Theta(\mu)$, according to Assumption 2, $\underline{\alpha}_V^{-1} \circ \tilde\gamma(|e|) \leq L_2 |e|$ where $\tilde\gamma(|e|) = \alpha^{-1}((1-\bar\sigma)^{-1} \gamma(|e|))$ here. Thus, we have that $L_2 |e| < |x|$ implies $\tilde\gamma(|e|) < \underline{\alpha}_V(|x|) \leq V(x)$ from Assumption 1, that ensures in return $\frac{\partial V}{\partial x} f(x, e) < -\bar\sigma \bar\alpha V(x)$ in view of (11) (with $\bar\sigma = (1-\sigma)$). Consequently, we see that after each jump, the time $L_2 |e|$ grows from 0 (since $e^+ = 0$) to $|x|$ ensures that $q \in C \backslash D_1$ since after each jump $V(x^+) = \eta^+ V(x^+ + e^+)$ and $V(x)$

will decrease faster than $\eta V(x+e)$ for some time. Hence, for all $(t_j, j) \preceq (t, j) \in \mathrm{dom}\,\phi$, $L_2 \frac{|e(t,j)|}{|x(t,j)|} < 1$ implies that $\phi(t,j) \in C \backslash D_1$. Noting that $\eta(t,j) > \varepsilon$ implies that $q \in C \backslash D_2$ for all $(t_j, j) \preceq (t, j) \in \mathrm{dom}\,\phi$, we deduce that item (iii-b) of Lemma 1 is satisfied with $b = 1$. We now prove that item (iii-c) of Lemma 1 holds by following similar lines as in (11) in [19]. The dynamics of $\frac{|e|}{|x|}$ are:

$$
\begin{aligned}
\frac{d}{dt} \frac{|e|}{|x|} &= \frac{d}{dt} \frac{(e^\mathrm{T} e)^{\frac{1}{2}}}{(x^\mathrm{T} x)^{\frac{1}{2}}} \\
&= \frac{(e^\mathrm{T} e)^{-\frac{1}{2}} e^\mathrm{T} \dot{e}(x^\mathrm{T} x)^{\frac{1}{2}} - (x^\mathrm{T} x)^{-\frac{1}{2}} x^\mathrm{T} \dot{x}(e^\mathrm{T} e)^{\frac{1}{2}}}{x^\mathrm{T} x} \\
&= \frac{e^\mathrm{T} \dot{e}}{|e||x|} - \frac{x^\mathrm{T} \dot{x}}{|x|^2} \frac{|e|}{|x|} \\
&\leq \frac{|e||\dot{e}|}{|e||x|} + \frac{|x||\dot{x}|}{|x|^2} \frac{|e|}{|x|} = \left(1 + \frac{|e|}{|x|}\right) \frac{|\dot{x}|}{|x|}
\end{aligned} \quad (37)
$$

using (14), $\frac{d}{dt} \frac{|e|}{|x|} \leq \left(1 + \frac{|e|}{|x|}\right) \frac{L_1|x| + L_1|e|}{|x|} = L_1 \left(1 + \frac{|e|}{|x|}\right)^2$. Therefore, noting that the derivative of $q \mapsto \frac{1-\eta}{1-\varepsilon}$ along the solutions to (20) is $q \mapsto \frac{1+\bar{\sigma}\bar{\alpha}}{1-\varepsilon}$, we deduce that item (iii-c) of Lemma 1 is verified with $\lambda : s \mapsto \max\{L_1 L_2(1+s^2), \frac{1+\bar{\sigma}\bar{\alpha}}{1-\varepsilon}\}$. Finally, by invoking Lemma 1, item (iv) of Theorem 1 holds, i.e. solutions to (20) have a semiglobal dwell time on $\mathbb{R}^{n_q} \backslash \mathcal{A}$, and we obtain that the set $\mathcal{A}$ is S-GAS using Theorem 1. $\square$

**Proof of Theorem 3.** The proof consists in checking the conditions of Theorem 1 for system (25) and then in applying it. We have that $G(D) = \{q \in \mathbb{R}^{n_q} : e = 0, \eta \geq 0\} \subset \mathbb{R}^{n_x+n_e} \times \mathbb{R}_{\geq 0} = C \cup D$. Consider system (25) and the following candidate Lyapunov function for $q \in \mathbb{R}^{n_q}$:

$$R(q) = \max\{V(x), W(e), \eta\}. \quad (38)$$

It can be verified that $R$ satisfies item (i) of Theorem 1 with $\upsilon(\eta) = \eta$, using Remark 2.3 in [11] and the fact that $C \cup D = \mathbb{R}^{n_x+n_e} \times \mathbb{R}_{\geq 0}$. On $C$, we have that $W(e) \leq \max\{V(x), \eta\} = R(q)$, therefore in view of (11) and since $\dot\eta = -\delta(\eta)$, we have that item (ii) of Theorem 1 is ensured with $\alpha_R(s) = \min\{(1-\sigma)\alpha(s), \delta(s)\}$ for $s \in \mathbb{R}_{\geq 0}$. Let $q \in D$, $R(G(q)) = \max\{V(x), W(0), W(e)\} = R(q)$: item (iii) of Theorem 1 is ensured. We now show that item (iv) of Theorem 1 holds using Lemma 1. We see that $G(0, e, \eta) = (0, 0, W(e)) \in \mathcal{A} \cup \mathcal{M}$ where $\mathcal{A} = \{q : (x, e, \eta) = 0\}$ and $\mathcal{M} = \{q \in C \cup D : (x, e) = 0 \text{ and } \eta \neq 0\}$. The set $\mathcal{M}$ is forward invariant for (25) as $f(0,0) = 0$ and $g(0,0) = 0$ in view of Assumption 2 and item (iii) of Theorem 3 and since $\eta(t,j) > 0$ for all $(t,j) \in \mathrm{dom}\,\phi$ for $\eta(0,0) > 0$. Let $\psi : q \mapsto \frac{|e|}{|x|}$ that is defined on $\Theta(\mu)$, where $\Theta(\mu)$ comes from Lemma 1 with $\mu \in \mathbb{R}_{>0}$. For any $q \in D$ such that $G(q) \in \Theta(\mu)$, $\psi(G(q)) = 0$ so item (iii-a) of Lemma 1 is satisfied with $a = 0$. By following similar lines than in the proof of Theorem 2, it can be shown that item (iii-b) of Lemma 1 is verified with $b = L_2^{-1}$. We now need to prove that item (iii-c) of Lemma 1 is guaranteed. We investigate the dynamics of $\frac{|e|}{|x|}$:

$$
\begin{aligned}
\frac{d}{dt} \frac{|e|}{|x|} &= \frac{d}{dt} \frac{(e^\mathrm{T} e)^{\frac{1}{2}}}{(x^\mathrm{T} x)^{\frac{1}{2}}} = \frac{(e^\mathrm{T} e)^{-\frac{1}{2}} e^\mathrm{T} \dot{e}(x^\mathrm{T} x)^{\frac{1}{2}} - (x^\mathrm{T} x)^{-\frac{1}{2}} x^\mathrm{T} \dot{x}(e^\mathrm{T} e)^{\frac{1}{2}}}{x^\mathrm{T} x} \\
&\leq \frac{|\dot{e}|}{|x|} + \frac{|\dot{x}||e|}{|x|^2},
\end{aligned} \quad (39)
$$

using (14) and item (iii) of Theorem 3, we obtain:

$$
\begin{aligned}
\frac{d}{dt} \frac{|e|}{|x|} &\leq \frac{L_3(|x|+|e|)}{|x|} + \frac{L_1(|x|+|e|)|e|}{|x|^2} \\
&= L_3 + L_3 \frac{|e|}{|x|} + L_1 \frac{|e|}{|x|} + L_1 \frac{|e|^2}{|x|^2},
\end{aligned} \quad (40)
$$

we see that item (iii-c) is guaranteed with $\lambda : s \mapsto L_3 + (L_1 + L_3)s + L_1 s^2$. By applying Lemma 1, item (iv) of Theorem 1 is verified: solutions to (25) have a semiglobal dwell time on $\mathbb{R}^{n_q} \backslash \{0\}$. As a consequence, the set $\mathcal{A}$ is S-GAS. $\square$

## References

[1] D. Angeli. A Lyapunov approach to incremental stability properties. *IEEE Transactions on Automatic Control*, 47(3):410–421, 2002.
[2] K.E. Arzén. A simple event-based PID controller. In $14^{th}$ *IFAC World Congress, Beijing, China*, 1999.
[3] K. J. Åström and B. Wittenmark. *Computer-controlled systems, theory and design*. Prentice-Hall, New Jersey, U.S.A., 3rd edition, 1996.
[4] K.J. Åström and B.M. Bernhardsson. Comparison of Riemann and Lebesgue sampling for first order stochastic systems. In *CDC (IEEE Conference on Decision and Control), Las Vegas, U.S.A.*, 2002.
[5] C. Cai, A.R. Teel, and R. Goebel. Smooth Lyapunov functions for hybrid systems part II: (pre)asymptotically stable compact sets. *IEEE Transactions on Automatic Control*, 53(3):734–748, 2008.
[6] M.C.F. Donkers and W.P.M.H. Heemels. Output-based event-triggered control with guaranteed $\mathcal{L}_\infty$-gain and improved event-triggering. In *IEEE Conference on Decision and Control, Atlanta, U.S.A.*, 2010.
[7] R. Goebel and A.R. Teel. Solution to hybrid inclusions via set and graphical convergence with stability theory applications. *Automatica*, 42:573–587, 2006.
[8] W.P.M.H. Heemels, J.H. Sandee, and P.P.J. van den Bosch. Analysis of event-driven controllers for linear systems. *International Journal of Control*, 81(4):571–590, 2009.
[9] H.K. Khalil. *Nonlinear Systems*. Prentice-Hall, Englewood Cliffs, New Jersey, U.S.A., 3rd edition, 2002.
[10] D. S. Laila and D. Nešić. Open and closed loop dissipation inequalities under sampling and controller emulation. *European Journal of Control*, 18:109–125, 2002.
[11] D.S. Laila and D. Nešić. Lyapunov based small-gain theorem for parameterized discrete-time interconnected ISS systems. In *CDC (IEEE Conference on Decision and Control), Las Vegas, U.S.A.*, pages 2292–2297, 2002.
[12] J. Löfberg. Pre- and post-processing sum-of-squares programs in practice. *IEEE Trans. on Automatic Control*, 54(5):1007–1011, 2009.
[13] D. Nešić, A.R. Teel, and D. Carnevale. Explicit computation of the sampling period in emulation of controllers for nonlinear sampled-data systems. *IEEE Trans. on Automatic Control*, 54(3):619–624, 2009.
[14] D. Nešić, L. Zaccarian, and A.R. Teel. Stability properties of reset systems. *Automatica*, 44:2019–2026, 2008.
[15] G. Otanez, J.R. Moyne, and D.M. Tilbury. Using deadbands to reduce communication in networked control systems. In *ACC (American Control Conference)*, 2002.
[16] R. Postoyan, A. Anta, D. Nešić, and P. Tabuada. A unifying Lyapunov-based framework for the event-triggered control of nonlinear systems. In *IEEE Conference on Decision and Control and European Control Conference, Orlando, U.S.A.*, 2011.
[17] R. Postoyan, P. Tabuada, D. Nešić, and A. Anta. Event-triggered and self-triggered stabilization of distributed networked control systems. In *CDC / ECC (IEEE Conference on Decision and Control and European Control Conference) Orlando, U.S.A.*, 2011.
[18] E.D. Sontag and Y. Wang. On characterizations of the input-to-state stability property. *Systems & Control Letters*, 24(5):351–359, 1995.
[19] P. Tabuada. Event-triggered real-time scheduling of stabilizing control tasks. *IEEE Trans. on Automatic Control*, 52(9):1680–1685, 2007.
[20] X. Wang and M.D. Lemmon. Event design in event-triggered feedback control systems. In *CDC (IEEE Conference on Decision and Control) Cancun, Mexico*, pages 2105–2110, 2008.